\documentclass[12pt]{article}
\usepackage{latexsym,amsfonts,amssymb}
\usepackage[dvips]{color}
\usepackage{colordvi,multicol}
\usepackage{amsmath}
\usepackage{graphicx}
\usepackage{pdflscape}
\usepackage{rotating}

\usepackage{latexsym,amsfonts,amssymb}
\usepackage[dvips]{color}
\usepackage{colordvi,multicol}
\usepackage{amsmath}
\usepackage{graphicx}
\usepackage{pdflscape}
\usepackage{rotating}
\usepackage{breqn}

\newtheorem{thm}{Theorem}[section]

\newtheorem{rem}[thm]{Remark}

\date{}
\begin{document}

\title{\bf Variation comparison between infinitely divisible distributions and the normal distribution}
\author{Ping Sun, Ze-Chun Hu  and Wei Sun\thanks{Corresponding author.}\\ \\
  {\small Business School, Chengdu University, Chengdu 610106, China}\\ \\
 {\small College of Mathematics, Sichuan University, Chengdu 610065, China}\\ \\
{\small Department of Mathematics and Statistics, Concordia University, Montreal H3G 1M8,  Canada}\\ \\
{\small  sunping@cdu.edu.cn\ \ \ \ zchu@scu.edu.cn\ \ \ \ wei.sun@concordia.ca}}

\maketitle

\begin{abstract}
\noindent
Let $X$ be a random variable with finite second moment. We investigate the inequality:  $P\{|X-E[X]|\le \sqrt{{\rm Var}(X)}\}\ge P\{|Z|\le 1\}$, where $Z$ is a standard normal random variable. We prove that this inequality holds for many familiar infinitely divisible continuous distributions including the Laplace, Gumbel, Logistic, Pareto, infinitely divisible Weibull,  log-normal, student's $t$ and  inverse Gaussian distributions. Numerical results are  given to show that the  inequality with continuity correction also holds for some infinitely divisible discrete distributions.
\end{abstract}



\section{Introduction}

Tomaszewski's conjecture says that if $T=\sum_{i=1}^na_ix_i$, where $\sum_{i=1}^n a_i^2=1$ and  $\{x_i\}$ is a sequence of independent $\{-1,1\}$-valued symmetric random variables, then $P\{|T|\leq 1\}\geq 1/2$. This conjecture has applications in probability theory, geometric analysis, computer science, economics and management science. Recently Keller and Klein \cite{KK22} completely solved Tomaszewski's conjecture. We refer the reader to Keller and  Klein \cite{KK22} for the details, and Dvorak and Klein \cite{DK22}  and Hu et al. \cite{Hu} for some related problems.
Motivated by Tomaszewski's conjecture,  we established the following result among other things in \cite{SHS}:
\begin{thm}\label{thm1100}(\cite{SHS})
Let $\alpha,\beta$ be arbitrary positive real numbers,  $X_{\alpha,\beta}$  be a Gamma random variable with shape parameter $\alpha$ and scale parameter $\beta$, and $Z$ be a standard normal random variable. Then,
\begin{equation}\label{inf111}
P\left\{|X_{\alpha,\beta}-E[X_{\alpha,\beta}]|\le \sqrt{{\rm Var}(X_{\alpha,\beta})}\right\}>P\{|Z|\le 1\}\approx 0.6827,
\end{equation}
and
$$
\inf_{\alpha,\beta}P\left\{|X_{\alpha,\beta}-E[X_{\alpha,\beta}]|\le \sqrt{{\rm Var}(X_{\alpha,\beta})}\right\}=P\{|Z|\le 1\}.
$$
\end{thm}

Note that the inequality of type (\ref{inf111}) does not hold for all continuous random variables. In particular, it does not hold for the uniform and the  Beta random variables. Let $a<b$ and  $X_{a,b}$ be a Uniform$(a,b)$  random variable. We have
\begin{eqnarray*}
P\left\{|X_{a,b}-E[X_{a,b}]|\le \sqrt{{\rm Var}(X_{a,b})}\right\}
&=&P\left\{\left|X_{a,b}-\frac{a+b}{2}\right|\le \sqrt{\frac{(b-a)^2}{12}}\right\}\\
&=&\frac{2}{\sqrt{12}}\\
&<&0.6827.
\end{eqnarray*}
 Let $\alpha,\beta>0$ and  $X_{\alpha,\beta}$ be a Beta random variable with parameters $\alpha$ and $\beta$. Define
$$
J_{\alpha,\beta}:=P\left\{|X_{\alpha,\beta}-E[X_{\alpha,\beta}]|\le \sqrt{{\rm Var}(X_{\alpha,\beta})}\right\}.
$$
Below is the graph of the function $(J_{2,\beta}-0.6827)$ for $\beta\in[1,20]$:

\begin{figure}[h]
\begin{center}
\scalebox{0.5}{\includegraphics{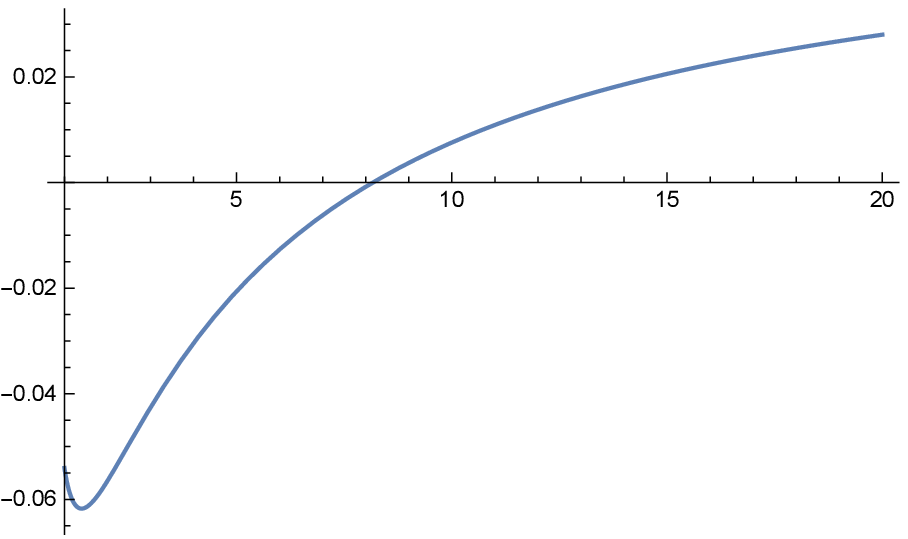}}
\end{center}
\end{figure}
\begin{center}
{\small Figure 1: Function $(J_{2,\beta}-0.6827)$ for $\beta\in [1,20]$.}
\end{center}

It is well-known that the Gamma distribution is infinitely divisible and, for any L\'evy process $\{L_t,t\ge0\}$, the distribution of $L_t$   is infinitely divisible. Inspired by Theorem \ref{thm1100}, it is natural to ask if any infinitely divisible  random variable $X$ with finite second moment satisfies the following inequality:
\begin{equation}\label{inf300}
P\left\{|X-E[X]|\le \sqrt{{\rm Var}(X)}\right\}\ge P\{|Z|\le 1\}.
\end{equation}
Recall that the empirical rule tells us in an approximately normal distribution about 68\% of the values fall within one standard deviation of the mean. It is very interesting and totally unexpected that actually  inequality (\ref{inf300}) can be established for many familiar infinitely divisible distributions.

In Sections 2--6, we will show that  (\ref{inf300}) holds for  the Laplace, Gumbel, Logistic, Pareto, infinitely divisible Weibull,  log-normal, student's $t$ and  inverse Gaussian distributions. We would like to point
out that, for the Weibull distribution,  (\ref{inf300}) might not holds if the distribution is not infinitely
divisible, i.e., if the parameter is bigger than 1 (cf. Remark \ref{May31}). Inequality  (\ref{inf300}) reveals a deep
relationship between some infinitely divisible distributions and the normal distribution. In Section 7, we consider (\ref{inf300}) for discrete random variables and make some remarks.  Although (\ref{inf300}) also holds for some infinitely divisible discrete distributions, e.g., the geometric distribution, it should be modified with continuity correction when general discrete distributions are considered.

\section{Laplace, Gumbel, Logistic and Pareto distributions}\setcounter{equation}{0}

\subsection{Laplace distribution}

\begin{thm}\label{thm25a}
Let $\mu\in\mathbb{R}$, $b>0$,  $X_{\mu,b}$ be a Laplace random variable with parameters $\mu$ and $b$, and $Z$ be a standard normal random variable. Then,
$$
P\left\{|X_{\mu,b}-E[X_{\mu,b}]|\le \sqrt{{\rm Var}(X_{\mu,b})}\right\}\approx0.7568833>P\{|Z|\le 1\}\approx 0.6827.
$$
\end{thm}

\noindent {\bf Proof.}\ \  We have (cf. \cite{WL})
$$
E[X_{\mu,b}]=\mu,\ \ \ \ {\rm Var}(X_{\mu,b})=2b^2,
$$
and
$$
P\{X_{\mu,b}\le x\}= \left\{ \begin{array}{ll}
         \frac{1}{2}\exp\left(\frac{x-\mu}{b}\right),\ \  & \mbox{if $x \le \mu$},\\
        1-\frac{1}{2}\exp\left(-\frac{x-\mu}{b}\right),\ \  & \mbox{if $x > \mu$}.\end{array} \right.
$$
Then,
\begin{eqnarray*}
P\left\{|X_{\mu,b}-E[X_{\mu,b}]|\le \sqrt{{\rm Var}(X_{\mu,b})}\right\}
&=&1-\frac{1}{2}\exp(-\sqrt{2})-\frac{1}{2}\exp(-\sqrt{2})\\
&\approx&0.7568833.
\end{eqnarray*}

\subsection{Gumbel distribution}

\begin{thm}\label{thm25b}
Let $\mu\in\mathbb{R}$, $\beta>0$,  $X_{\mu,\beta}$ be a Gumbel random variable with parameters $\mu$ and $\beta$, and $Z$ be a standard normal random variable. Then,
$$
P\left\{|X_{\mu,\beta}-E[X_{\mu,\beta}]|\le \sqrt{{\rm Var}(X_{\mu,\beta})}\right\}\approx0.723751>P\{|Z|\le 1\}\approx 0.6827.
$$
\end{thm}

\noindent {\bf Proof.}\ \  We have (cf. \cite{WG})
$$
E[X_{\mu,\beta}]=\mu+\beta\gamma,\ \ \ \ {\rm Var}(X_{\mu,\beta})=\frac{\pi^2\beta^2}{6},
$$
where $\gamma$ is the Euler constant, and
$$
P\{X_{\mu,\beta}\le x\}=e^{-e^{-(x-\mu)/\beta}}.
$$ Then,
\begin{eqnarray*}
P\left\{|X_{\mu,b}-E[X_{\mu,b}]|\le \sqrt{{\rm Var}(X_{\mu,b})}\right\}
&=&e^{-e^{-(\gamma+\pi/\sqrt{6})}}-e^{-e^{-(\gamma-\pi/\sqrt{6})}}\\
&\approx&0.723751.
\end{eqnarray*}

\subsection{Logistic distribution}

\begin{thm}\label{thm25c}
Let $\mu\in\mathbb{R}$, $s>0$,  $X_{\mu,s}$ be a Logistic random variable with parameters $\mu$ and $s$, and $Z$ be a standard normal random variable. Then,
$$
P\left\{|X_{\mu,s}-E[X_{\mu,s}]|\le \sqrt{{\rm Var}(X_{\mu,s})}\right\}\approx0.719641>P\{|Z|\le 1\}\approx 0.6827.
$$
\end{thm}

\noindent {\bf Proof.}\ \  We have (cf. \cite{WLo})
$$
E[X_{\mu,s}]=\mu,\ \ \ \ {\rm Var}(X_{\mu,s})=\frac{\pi^2s^2}{3},
$$
and
$$
P\{X_{\mu,s}\le x\}=\frac{1}{1+e^{-(x-\mu)/s}}.
$$ Then,
\begin{eqnarray*}
P\left\{|X_{\mu,s}-E[X_{\mu,s}]|\le \sqrt{{\rm Var}(X_{\mu,s})}\right\}
&=&\frac{1}{1+e^{-\pi/\sqrt{3}}}-\frac{1}{1+e^{\pi/\sqrt{3}}}\\
&\approx&0.719641.
\end{eqnarray*}

\subsection{Pareto distribution}

\begin{thm}\label{thm25d}
Let $x_m>0$, $\alpha>2$,  $X_{x_m,\alpha}$ be a Pareto random variable with parameters $x_m$ and $\alpha$, and $Z$ be a standard normal random variable. Then, for any $x_m>0$,
\begin{eqnarray*}
P\left\{|X_{x_m,\alpha}-E[X_{x_m,\alpha}]|\le \sqrt{{\rm Var}(X_{x_m,\alpha})}\right\}&\downarrow& 1-e^{-2}\quad {\rm as}\quad \alpha\uparrow\infty\\
&\approx&  0.8646647\\
&>&P\{|Z|\le 1\}\approx 0.6827.
\end{eqnarray*}
\end{thm}

\noindent {\bf Proof.}\ \  We have (cf. \cite{WP})
$$
E[X_{x_m,\alpha}]=\frac{\alpha x_m}{\alpha-1},\ \ \ \ {\rm Var}(X_{x_m,\alpha})=\frac{\alpha x^2_m}{(\alpha-1)^2(\alpha-2)},
$$
and
$$
P\{X_{x_m,\alpha}\le x\}=1-\left(\frac{x_m}{x}\right)^{\alpha},\ \ \ \ x\ge x_m.
$$ Then,
\begin{eqnarray*}
&&P\left\{|X_{x_m,\alpha}-E[X_{x_m,\alpha}]|\le \sqrt{{\rm Var}(X_{x_m,\alpha})}\right\}\\
&=&1-\left(\frac{1}{\frac{\alpha}{\alpha-1}+\frac{\sqrt{\alpha}}{(\alpha-1)\sqrt{\alpha-2}}}\right)^{\alpha}\\
&=&1-\left[1+\frac{1}{(\alpha-1)(1-[1+(1-\frac{2}{\alpha})^{\frac{1}{2}}]^{-1})}\right]^{-\alpha}\\
&\rightarrow&1-e^{-2}\ \ {\rm as}\ \alpha\rightarrow\infty.
\end{eqnarray*}
Therefore, the proof is complete by noting that
$$
1+\frac{1}{(\alpha-1)(1-[1+(1-\frac{2}{\alpha})^{\frac{1}{2}}]^{-1})}
$$
is a strictly decreasing function of $\alpha\in(2,\infty)$.

\section{Weibull distribution}\setcounter{equation}{0}

\begin{thm}\label{thm25e}
Let $\lambda>0$, $0<k\le 1$,  $X_{\lambda,k}$ be a Weibull  random variable with parameters $\lambda$ and $k$, and $Z$ be a standard normal random variable. Then,
\begin{eqnarray*}
P\left\{|X_{\lambda,k}-E[X_{\lambda,k}]|\le \sqrt{{\rm Var}(X_{\lambda,k})}\right\}>P\{|Z|\le 1\}\approx 0.6827.
\end{eqnarray*}
\end{thm}

\noindent {\bf Proof.}\ \  We have (cf. \cite{WW})
$$
E[X_{\lambda,k}]=\lambda\Gamma\left(1+\frac{1}{k}\right),\ \ \ \ {\rm Var}(X_{\lambda,k})=\lambda^2\left[\Gamma\left(1+\frac{2}{k}\right)-\left\{\Gamma\left(1+\frac{1}{k}\right)\right\}^2\right],
$$
and
$$
P\{X_{\lambda,k}\le x\}=1-e^{-(x/\lambda)^k},\ \ \ \ x\ge 0.
$$ Then,
\begin{eqnarray*}
&&P\left\{|X_{\lambda,k}-E[X_{\lambda,k}]|\le \sqrt{{\rm Var}(X_{\lambda,k})}\right\}\\
&=&e^{-\left[\max\left\{0,\,\Gamma\left(1+\frac{1}{k}\right)-
\sqrt{\Gamma\left(1+\frac{2}{k}\right)-
\left\{\Gamma\left(1+\frac{1}{k}\right)\right\}^2}\right\}\right]^k}-e^{-\left[\Gamma\left(1+\frac{1}{k}\right)+\sqrt{\Gamma\left(1+\frac{2}{k}\right)-\left\{\Gamma\left(1+\frac{1}{k}\right)\right\}^2}\right]^k}\\
&:=&W_k.
\end{eqnarray*}

By the Legendre duplication formula, we obtain that for $k\in(0,1]$,
\begin{eqnarray*}
\frac{2\{\Gamma\left(1+\frac{1}{k}\right)\}^2}{\Gamma\left(1+\frac{2}{k}\right)}&=&\frac{2\{\frac{1}{k}\Gamma\left(\frac{1}{k}\right)\}^2}{\frac{2}{k}\Gamma\left(\frac{2}{k}\right)}\nonumber\\
&=&\frac{1}{k2^{\frac{2}{k}-1}}B\left(\frac{1}{2},\frac{1}{k}\right)\\
&=&\frac{1}{k2^{\frac{2}{k}-1}}\int_0^1t^{-\frac{1}{2}}(1-t)^{\frac{1}{k}-1}dt\\
&<&\frac{1}{k2^{\frac{2}{k}-1}}\int_0^1t^{-\frac{1}{2}}dt\\
&=&\frac{1}{k2^{\frac{2}{k}-2}}\\
&\le&1.
\end{eqnarray*}
Then, $\Gamma\left(1+\frac{1}{k}\right)-
\sqrt{\Gamma\left(1+\frac{2}{k}\right)-
\left\{\Gamma\left(1+\frac{1}{k}\right)\right\}^2}< 0$. Thus,
\begin{eqnarray*}
W_k&=&1-e^{-\left[\Gamma\left(1+\frac{1}{k}\right)+\sqrt{\Gamma\left(1+\frac{2}{k}\right)-\left\{\Gamma\left(1+\frac{1}{k}\right)\right\}^2}\right]^k}\\
&\ge&1-e^{-\left[2\Gamma\left(1+\frac{1}{k}\right)\right]^k}\\
&>&1-e^{-\left(2n!\right)^{\frac{1}{n+1}}}\ \ \ \ {\rm if}\ k\in\left(\frac{1}{n+1},\frac{1}{n}\right]\\
&\ge&1-e^{-\sqrt{2}}\\
&\approx&0.7568833,
\end{eqnarray*}
where the last inequality holds since
$$
\frac{\left(2n!\right)^{\frac{1}{n+1}}}{\left[2(n+1)!\right]^{\frac{1}{n+2}}}=\left(\frac{2n!}{(n+1)^{n+1}}\right)^{\frac{1}{(n+1)(n+2)}}<1,\ \ \ \ \forall n\in\mathbb{N}.
$$

\begin{rem}\label{May31}
It is known that the Weibull distribution is infinitely divisible if and only if $k\in(0,1]$ (cf. \cite[Remark 8.12, page 46 and E.29.10, page 194]{S}). In general, the following variation comparison inequality
$$
P\left\{|X_{\lambda,k}-E[X_{\lambda,k}]|\le \sqrt{{\rm Var}(X_{\lambda,k})}\right\}\ge P\{|Z|\le 1\}\approx 0.6827
$$
does not hold if $k>1$. For example, $W_3=0.667713<0.6827$. Below is the graph of the function $(W_k-0.6827)$ for $k\in [1,10]$.

\newpage\begin{figure}[h]
\begin{center}
\scalebox{0.5}{\includegraphics{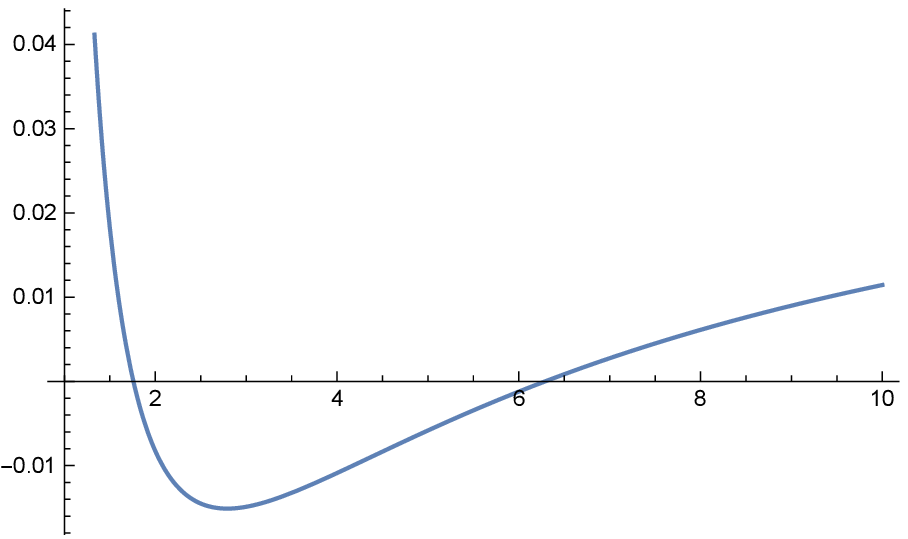}}
\end{center}
\end{figure}
\begin{center}
{\small Figure 2: Function $(W_k-0.6827)$ for $k\in [1,10]$.}
\end{center}
\end{rem}

\section{Log-normal distribution}\setcounter{equation}{0}

A log-normal  distribution is the probability distribution of a continuous random variable whose logarithm is normally distributed. Let  $Z$ be a standard normal variable, $\mu\in\mathbb{R}$ and $\sigma>0$. Then, the distribution of the random variable
$$X_{\mu,\sigma}=e^{\mu +\sigma Z}
$$
is called the log-normal distribution with parameters $\mu$ and $\sigma$. It is known that (cf. \cite{Wiki}) the density function of $X_{\mu,\sigma}$ is given by
$$
f_{\mu,\sigma}(x)=\frac{1}{\sqrt{2\pi}\sigma x}\exp\left(-\frac{(\ln x-\mu)^2}{2\sigma^2}\right),\ \ \ \ x>0,
$$
and
$$
E[X_{\mu,\sigma}]=\exp\left(\mu+\frac{\sigma^2}{2}\right),\ \ \ \ {\rm Var}(X_{\mu,\sigma})=[\exp(\sigma^2)-1]\exp\left(2\mu+\sigma^2\right).
$$

\begin{thm}\label{thm11}
Let $\mu\in\mathbb{R}$, $\sigma>0$,  $X_{\mu,\sigma}$ be a log-normal random variable with parameters $\mu$ and $\sigma$, and $Z$ be a standard normal random variable. Then,
$$
P\left\{|X_{\mu,\sigma}-E[X_{\mu,\sigma}]|\le \sqrt{{\rm Var}(X_{\mu,\sigma})}\right\}>P\{|Z|\le 1\}\approx 0.6827,
$$
and
$$
\inf_{\mu,\sigma}P\left\{|X_{\mu,\sigma}-E[X_{\mu,\sigma}]|\le \sqrt{{\rm Var}(X_{\mu,\sigma})}\right\}=P\{|Z|\le 1\}.
$$
\end{thm}

\noindent {\bf Proof.}\ \ Denote by $\Phi$ the cumulative distribution function of the standard normal distribution. We have
\begin{eqnarray*}
&&P\left\{|X_{\mu,\sigma}-E[X_{\mu,\sigma}]|\le \sqrt{{\rm Var}(X_{\mu,\sigma})}\right\}\\
&=&P\left\{\left|X_{\mu,\sigma}-\exp\left(\mu+\frac{\sigma^2}{2}\right)\right|\le \sqrt{\exp(\sigma^2)-1}\exp\left(\mu+\frac{\sigma^2}{2}\right)\right\}\\
&=&\left\{\begin{array}{ll}\Phi\left(\frac{\frac{\sigma^2}{2}+\ln(1+\sqrt{\exp(\sigma^2)-1})}{\sigma}\right),\ \ \ \  &\sigma\ge\sqrt{\ln 2}, \\ \Phi\left(\frac{\frac{\sigma^2}{2}+\ln(1+\sqrt{\exp(\sigma^2)-1})}{\sigma}\right)-\Phi\left(\frac{\frac{\sigma^2}{2}+\ln(1-\sqrt{\exp(\sigma^2)-1})}{\sigma}\right),\ \ \ \  &
0<\sigma<\sqrt{\ln 2}.\end{array}\right.
\end{eqnarray*}
By L'H$\hat{o}$pital's rule, we can show that
\begin{eqnarray*}
&&\lim_{\sigma\rightarrow0}\Bigg[\Phi\left(\frac{\frac{\sigma^2}{2}+\ln(1+\sqrt{\exp(\sigma^2)-1})}{\sigma}\right)-\Phi\left(\frac{\frac{\sigma^2}{2}+\ln(1-\sqrt{\exp(\sigma^2)-1})}{\sigma}\right)\Bigg]\\
&=&\Phi(1)-\Phi(-1)\\
&\approx &0.6827.
\end{eqnarray*}
Hence, to complete the proof, we need only show that the function
$$
\frac{\frac{\sigma^2}{2}+\ln(1+\sqrt{\exp(\sigma^2)-1})}{\sigma}
$$
is increasing on $\sigma\in(0,\infty)$ and the function
$$
\frac{\frac{\sigma^2}{2}+\ln(1-\sqrt{\exp(\sigma^2)-1})}{\sigma}
$$
is decreasing on $\sigma\in(0,\sqrt{\ln 2})$.

Define
$$
y=\sqrt{\exp(\sigma^2)-1}.
$$
Then, for $0<\sigma<\infty\Leftrightarrow0<y<\infty$,
$$
\frac{\frac{\sigma^2}{2}+\ln(1+\sqrt{\exp(\sigma^2)-1})}{\sigma}=\frac{\frac{\ln(1+y^2)}{2}+\ln(1+y)}{\sqrt{\ln(1+y^2)}}=\frac{\ln\{\sqrt{1+y^2}\cdot(1+y)\}}{\sqrt{\ln(1+y^2)}},
$$
and  for $0<\sigma<\sqrt{\ln 2}\Leftrightarrow0<y<1$,
$$
\frac{\frac{\sigma^2}{2}+\ln(1-\sqrt{\exp(\sigma^2)-1})}{\sigma}=\frac{\frac{\ln(1+y^2)}{2}+\ln(1-y)}{\sqrt{\ln(1+y^2)}}=\frac{\ln\{\sqrt{1+y^2}\cdot(1-y)\}}{\sqrt{\ln(1+y^2)}}.
$$
We have
\begin{eqnarray*}
\Bigg[\frac{\ln\{\sqrt{1+y^2}(1+y)\}}{\sqrt{\ln(1+y^2)}}\Bigg]'=\frac{\ln[(1+y^2)^{\frac{3y^2+y+2}{2}}(1+y)^{-(y+y^2)}]}{(1+y)(1+y^2)[\ln(1+y^2)]^{3/2}},
\end{eqnarray*}
and
\begin{eqnarray*}
-\Bigg[\frac{\ln\{\sqrt{1+y^2}(1-y)\}}{\sqrt{\ln(1+y^2)}}\Bigg]'=\frac{-\ln[(1+y^2)^{\frac{-3y^2+y-2}{2}}(1-y)^{y^2-y}]}{(1-y)(1+y^2)[\ln(1+y^2)]^{3/2}}.
\end{eqnarray*}

We will show that the above two derivatives are positive. We have
\begin{eqnarray*}
\ln[(1+y^2)^{\frac{3y^2+y+2}{2}}(1+y)^{-(y+y^2)}]=\ln(1+y^2)\cdot\frac{3y^2+y+2}{2}-y(1+y)\cdot\ln(1+y),
\end{eqnarray*}
which is obviously positive if $y\ge 1$; and if $0<y<1$,
\begin{eqnarray*}
&&\ln(1+y^2)\cdot\frac{3y^2+y+2}{2}-y(1+y)\cdot\ln(1+y)\\
&>&\left(y^2-\frac{y^4}{2}+\frac{y^6}{3}-\frac{y^8}{4}\right)\frac{3y^2+y+2}{2}-y(1+y)\left(y-\frac{y^2}{2}+\frac{y^3}{3}\right)\\
&=&\frac{y^4}{24}(28 - 14 y - 10 y^2 + 4 y^3 + 6 y^4 - 3 y^5 - 9 y^6)\\
&>&\frac{y^4}{24}(28-14-10-3)\\
&>&0.
\end{eqnarray*}
For  $0<y<1$, we have
\begin{eqnarray*}
&&-\ln[(1+y^2)^{\frac{-3y^2+y-2}{2}}(1-y)^{y^2-y}]\\
&=&\ln(1+y^2)\cdot{\frac{3y^2-y+2}{2}}+y(1-y)\cdot\ln(1-y)\\
&>&\left(y^2-\frac{y^4}{2}\right)\frac{3y^2-y+2}{2}+y\left(-y+\frac{y^2}{2}\right)\\
&=&\frac{y^4}{4}(4 + y - 3 y^2)\\
&>&0.
\end{eqnarray*}
Therefore, the proof is complete.

\section{Student's $t$-distribution}\setcounter{equation}{0}

Let $3\le \nu\in\mathbb{N}$ and $X_{\nu}$ be a $t$-random variable with $\nu$   degrees of freedom. Denote by $F(a, b; c; z)$ the hypergeometric function (cf. \cite{R}):
$$
    F(a, b; c; z)=\sum_{j=0}^{\infty}\frac{(a)_{j}(b)_{j}}{(c)_{j}}\cdot\frac{z^j}{j!},\ \ \ \ |z|<1,
$$
where $(\alpha)_{j}:=\alpha(\alpha+1)\cdots(\alpha+j-1)$ for $j\ge 1$, and $(\alpha)_0=1$ for $\alpha\not=0$. The density function and the cumulative distribution function of $X_{\nu}$  are given by (cf. \cite{Wiki3})
\begin{equation}\label{eqn12}
f_{\nu}(x)=\frac{\Gamma(\frac{\nu+1}{2})}{\sqrt{\nu\pi}\Gamma(\frac{\nu}{2})}\left(1+\frac{x^2}{\nu}\right)^{-\frac{\nu+1}{2}},\ \ \ \ x\in\mathbb{R},
\end{equation}
\begin{equation}\label{eqn13}
F_{\nu}(x)=\frac{1}{2}+x\Gamma\left(\frac{\nu+1}{2}\right)\frac{F\left(\frac{1}{2},\frac{\nu+1}{2};\frac{3}{2};-\frac{x^2}{\nu}\right)}{\sqrt{\nu\pi}\Gamma(\frac{\nu}{2})},\ \ \ \ x\in\mathbb{R},
\end{equation}
respectively, and
$$
E[X_{\nu}]=0,\ \ \ \ {\rm Var}(X_{\nu})=\frac{\nu}{\nu-2}.
$$

\begin{thm}
Let $\nu\ge3$,  $X_{\nu}$ be a $t$-random variable with $\nu$   degrees of freedom, and $Z$ be a standard normal random variable. Then,
$$
P\left\{|X_{\nu}-E[X_{\nu}]|\le \sqrt{{\rm Var}(X_{\nu})}\right\}>P\{|Z|\le 1\}\approx 0.6827,
$$
and
$$
\inf_{\nu}P\left\{|X_{\nu}-E[X_{\nu}]|\le \sqrt{{\rm Var}(X_{\nu})}\right\}=P\{|Z|\le 1\}.
$$
\end{thm}

\noindent {\bf Proof.}\ \ By (\ref{eqn13}), we get
\begin{eqnarray*}
P\left\{|X_{\nu}-E[X_{\nu}]|\le \sqrt{{\rm Var}(X_{\nu})}\right\}
&=&2\sqrt{\frac{\nu}{\nu-2}}\frac{\Gamma(\frac{\nu+1}{2})}{\sqrt{\nu\pi}\Gamma(\frac{\nu}{2})}F\left(\frac{1}{2},\frac{\nu+1}{2};\frac{3}{2};-\frac{1}{\nu-2}\right)\\
&:=&J_{\nu}.
\end{eqnarray*}
It is well-known  that student's $t$-distribution converges to the standard normal distribution as $\nu$ tends to infinity (cf. \cite[Page 453]{W}). Then,
$$
\lim_{\nu\rightarrow\infty}P\left\{|X_{\nu}-E[X_{\nu}]|\le \sqrt{{\rm Var}(X_{\nu})}\right\}=P\{|Z|\le 1\}.
$$
To complete the proof, we need only show that
$$
\frac{J_{\nu+2}}{J_{\nu}}<1,\ \ \ \ \forall\nu\ge 3.
$$

We have
\begin{eqnarray}\label{19A}
\frac{J_{\nu+2}}{J_{\nu}}<1&
\Leftrightarrow&\frac{(\nu+1)(\nu-2)^{\frac{1}{2}}}{\nu^{\frac{3}{2}}}\frac{F\left(\frac{1}{2},\frac{\nu+3}{2};\frac{3}{2};-\frac{1}{\nu}\right)}{F\left(\frac{1}{2},\frac{\nu+1}{2};\frac{3}{2};-\frac{1}{\nu-2}\right)}<1\nonumber\\
&
\Leftrightarrow&F\left(\frac{1}{2},\frac{\nu+3}{2};\frac{3}{2};-\frac{1}{\nu}\right)<\frac{\nu^{\frac{3}{2}}}{(\nu+1)(\nu-2)^{\frac{1}{2}}}F\left(\frac{1}{2},\frac{\nu+1}{2};\frac{3}{2};-\frac{1}{\nu-2}\right).
\end{eqnarray}
By the relation of Gauss between contiguous functions (cf. \cite[Page 71, 21(13)]{R}) and the fact that $F(a,b;a;z)=(1-z)^{-b}$, we get
\begin{eqnarray*}
\frac{\nu+1}{2}F\left(\frac{1}{2},\frac{\nu+3}{2};\frac{3}{2};-\frac{1}{\nu}\right)
&=&\frac{\nu}{2}F\left(\frac{1}{2},\frac{\nu+1}{2};\frac{3}{2};-\frac{1}{\nu}\right)+\frac{1}{2}F\left(\frac{1}{2},\frac{\nu+1}{2};\frac{1}{2};-\frac{1}{\nu}\right)\\
&=&\frac{\nu}{2}F\left(\frac{1}{2},\frac{\nu+1}{2};\frac{3}{2};-\frac{1}{\nu}\right)+\frac{1}{2}\left(\frac{\nu}{\nu+1}\right)^{\frac{\nu+1}{2}},
\end{eqnarray*}
which implies that
\begin{eqnarray*}
F\left(\frac{1}{2},\frac{\nu+3}{2};\frac{3}{2};-\frac{1}{\nu}\right)=\frac{\nu}{\nu+1}F\left(\frac{1}{2},\frac{\nu+1}{2};\frac{3}{2};-\frac{1}{\nu}\right)+\frac{1}{\nu+1}\left(\frac{\nu}{\nu+1}\right)^{\frac{\nu+1}{2}}.
\end{eqnarray*}
Then, by (\ref{eqn12})--(\ref{19A}), we get
\begin{eqnarray}\label{19bn}
&&\frac{J_{\nu+2}}{J_{\nu}}<1\nonumber\\
&\Leftrightarrow&\ \ \frac{\nu}{\nu+1}F\left(\frac{1}{2},\frac{\nu+1}{2};\frac{3}{2};-\frac{1}{\nu}\right)+\frac{1}{\nu+1}\left(\frac{\nu}{\nu+1}\right)^{\frac{\nu+1}{2}}<\frac{\nu^{\frac{3}{2}}}{(\nu+1)(\nu-2)^{\frac{1}{2}}}F\left(\frac{1}{2},\frac{\nu+1}{2};\frac{3}{2};-\frac{1}{\nu-2}\right)\nonumber\\
&\Leftrightarrow&\ \ F\left(\frac{1}{2},\frac{\nu+1}{2};\frac{3}{2};-\frac{1}{\nu}\right)+\frac{1}{\nu}\left(\frac{\nu}{\nu+1}\right)^{\frac{\nu+1}{2}}<\left(\frac{\nu}{\nu-2}\right)^{\frac{1}{2}}F\left(\frac{1}{2},\frac{\nu+1}{2};\frac{3}{2};-\frac{1}{\nu-2}\right)\nonumber\\
&\Leftrightarrow&\ \ \int_0^1f_{\nu}(x)dx+\frac{\Gamma(\frac{\nu+1}{2})}{\sqrt{\nu\pi}\Gamma(\frac{\nu}{2})}\cdot\frac{1}{\nu}\left(\frac{\nu}{\nu+1}\right)^{\frac{\nu+1}{2}}<\int_0^{\left(\frac{\nu}{\nu-2}\right)^{\frac{1}{2}}}f_{\nu}(x)dx\nonumber\\
&\Leftrightarrow&\frac{\Gamma(\frac{\nu+1}{2})}{\sqrt{\nu\pi}\Gamma(\frac{\nu}{2})}\cdot\frac{1}{\nu}\left(\frac{\nu}{\nu+1}\right)^{\frac{\nu+1}{2}}<\int_1^{\left(\frac{\nu}{\nu-2}\right)^{\frac{1}{2}}}\frac{\Gamma(\frac{\nu+1}{2})}{\sqrt{\nu\pi}\Gamma(\frac{\nu}{2})}\left(1+\frac{x^2}{\nu}\right)^{-\frac{\nu+1}{2}}dx\nonumber\\
&\Leftrightarrow&\frac{1}{\nu}\left(\frac{\nu}{\nu+1}\right)^{\frac{\nu+1}{2}}<\int_1^{\left(\frac{\nu}{\nu-2}\right)^{\frac{1}{2}}}\left(1+\frac{x^2}{\nu}\right)^{-\frac{\nu+1}{2}}dx\nonumber\\
&\Leftrightarrow&1<\nu\int_1^{\left(\frac{\nu}{\nu-2}\right)^{\frac{1}{2}}}\left(\frac{\nu+1}{\nu+x^2}\right)^{\frac{\nu+1}{2}}dx.
\end{eqnarray}

Below we show that inequality (\ref{19bn}) holds.
For $\nu=3$, we have
\begin{eqnarray*}
\nu\int_1^{\left(\frac{\nu}{\nu-2}\right)^{\frac{1}{2}}}\left(\frac{\nu+1}{\nu+x^2}\right)^{\frac{\nu+1}{2}}dx
&=&3\int_1^{3^{\frac{1}{2}}}\left(\frac{4}{3+x^2}\right)^{2}dx\\
&>&3\left\{(1.5-1)\left(\frac{4}{3+1.5^2}\right)^{2}+(3^{\frac{1}{2}}-1.5)\left(\frac{4}{3+3}\right)^{2}\right\}\\
&=& 1.180149\\
&>&1.
\end{eqnarray*}
For $\nu\ge 4$, we have
\begin{eqnarray*}
1<\nu\int_1^{\left(\frac{\nu}{\nu-2}\right)^{\frac{1}{2}}}\left(\frac{\nu+1}{\nu+x^2}\right)^{\frac{\nu+1}{2}}dx\Leftarrow1<\nu\left[\left(\frac{\nu}{\nu-2}\right)^{\frac{1}{2}}-1\right]\left[\frac{(\nu+1)(\nu-2)}{\nu(\nu-1)}\right]^{\frac{\nu+1}{2}}.
\end{eqnarray*}
Note that
\begin{eqnarray*}
&&4\left[\left(\frac{4}{4-2}\right)^{\frac{1}{2}}-1\right]\left[\frac{(4+1)(4-2)}{4(4-1)}\right]^{\frac{4+1}{2}}=1.050343>1,\\
&&5\left[\left(\frac{5}{5-2}\right)^{\frac{1}{2}}-1\right]\left[\frac{(5+1)(5-2)}{5(5-1)}\right]^{\frac{5+1}{2}}=1.060675>1,
\end{eqnarray*}
and
\begin{eqnarray*}
&&1<\nu\left[\left(\frac{\nu}{\nu-2}\right)^{\frac{1}{2}}-1\right]\left[\frac{(\nu+1)(\nu-2)}{\nu(\nu-1)}\right]^{\frac{\nu+1}{2}}\\
&\Leftrightarrow&\nu\left[\left(1+\frac{2}{\nu-2}\right)^{\frac{1}{2}}-1\right]\left[1-\frac{2}{\nu(\nu-1)}\right]^{\frac{\nu+1}{2}}-1>0\\
&\Leftarrow&\nu\left[\frac{1}{\nu-2}-\frac{1}{2(\nu-2)^2}\right]\left[1-\frac{\nu+1}{\nu(\nu-1)}\right]-1>0\\
&\Leftrightarrow&\frac{(2\nu-5)(\nu^2-2\nu-1)}{2(\nu-2)^2(\nu-1)}-1>0\\
&\Leftrightarrow&\nu^2-8\nu+13>0\\
&\Leftrightarrow&(\nu-4)^2-3>0,
\end{eqnarray*}
which obviously holds for  $\nu\ge 6$. Therefore, the proof is complete.

\section{Inverse Gaussian distribution}\setcounter{equation}{0}

Let $\mu,\lambda>0$ and $X_{\mu,\lambda}$ be an inverse Gaussian (also known as Wald) random variable with mean $\mu$ and shape parameter $\lambda$. The density function of $X_{\mu,\lambda}$  is given by (cf. \cite{Wiki2})
$$
f_{\mu,\lambda}(x)=\sqrt{\frac{\lambda}{2\pi x^3}}\exp\left(-\frac{\lambda(x-\mu)^2}{2\mu^2 x}\right),\ \ \ \ x>0,
$$
and
$$
E[X_{\mu,\lambda}]=\mu,\ \ \ \ {\rm Var}(X_{\mu,\lambda})=\frac{\mu^3}{\lambda}.
$$

\begin{thm}\label{thm11}
Let $\mu,\lambda>0$,  $X_{\mu,\lambda}$ be an inverse Gaussian random variable with parameters $\mu$ and $\lambda$, and $Z$ be a standard normal random variable. Then,
$$
P\left\{|X_{\mu,\lambda}-E[X_{\mu,\lambda}]|\le \sqrt{{\rm Var}(X_{\mu,\lambda})}\right\}>P\{|Z|\le 1\}\approx 0.6827,
$$
and
$$
\inf_{\mu,\lambda}P\left\{|X_{\mu,\lambda}-E[X_{\mu,\lambda}]|\le \sqrt{{\rm Var}(X_{\mu,\lambda}}\right\}=P\{|Z|\le 1\}.
$$
\end{thm}

\noindent {\bf Proof.}\ \ Denote by $\Phi$ the cumulative distribution function of the standard normal distribution. We have (cf. \cite{Wiki2})
\begin{eqnarray*}
&&P\left\{|X_{\mu,\lambda}-E[X_{\mu,\lambda}]|\le \sqrt{{\rm Var}(X_{\mu,\lambda})}\right\}\\
&=&P\left\{\left|X_{\mu,\lambda}-\mu\right|\le \sqrt{\frac{\mu^3}{\lambda}}\right\}\\
&=&\left\{\begin{array}{ll} \Phi\left(\frac{1}{(1+(\frac{\mu}{\lambda})^{\frac{1}{2}})^{\frac{1}{2}}}\right)+\exp\left(\frac{2\lambda}{\mu}\right)\Phi\left(-\frac{1+2(\frac{\lambda}{\mu})^{\frac{1}{2}}}{(1+(\frac{\mu}{\lambda})^{\frac{1}{2}})^{\frac{1}{2}}}\right),\ \ \ \  &\mu\ge\lambda, \\ \Phi\left(\frac{1}{(1+(\frac{\mu}{\lambda})^{\frac{1}{2}})^{\frac{1}{2}}}\right)+\exp\left(\frac{2\lambda}{\mu}\right)\Phi\left(-\frac{1+2(\frac{\lambda}{\mu})^{\frac{1}{2}}}{(1+(\frac{\mu}{\lambda})^{\frac{1}{2}})^{\frac{1}{2}}}\right) &\\
-\Phi\left(-\frac{1}{(1-(\frac{\mu}{\lambda})^{\frac{1}{2}})^{\frac{1}{2}}}\right)-\exp\left(\frac{2\lambda}{\mu}\right)\Phi\left(\frac{1-2(\frac{\lambda}{\mu})^{\frac{1}{2}}}{(1-(\frac{\mu}{\lambda})^{\frac{1}{2}})^{\frac{1}{2}}}\right),\ \ \ \  &
\mu<\lambda.\end{array}\right.
\end{eqnarray*}
Note that
\begin{eqnarray*}
&&\exp\left(\frac{2\lambda}{\mu}\right)\Phi\left(-\frac{1+2(\frac{\lambda}{\mu})^{\frac{1}{2}}}{(1+(\frac{\mu}{\lambda})^{\frac{1}{2}})^{\frac{1}{2}}}\right)\\
&=&\frac{1}{\sqrt{2\pi}}\exp\left(\frac{2\lambda}{\mu}\right)\int_{\frac{1+2(\frac{\lambda}{\mu})^{\frac{1}{2}}}{(1+(\frac{\mu}{\lambda})^{\frac{1}{2}})^{\frac{1}{2}}}}^{\infty}\exp\left(-\frac{y^2}{2}\right)dy\\
&=&\frac{1}{\sqrt{2\pi}}\exp\left(\frac{2\lambda}{\mu}\right)\int_{0}^{\infty}\exp\left(-\frac{\left(y+\frac{1+2(\frac{\lambda}{\mu})^{\frac{1}{2}}}{(1+(\frac{\mu}{\lambda})^{\frac{1}{2}})^{\frac{1}{2}}}\right)^2}{2}\right)dy\\
&=&\frac{1}{\sqrt{2\pi}}\exp\left(-\frac{1}{2[1+(\frac{\mu}{\lambda})^{\frac{1}{2}}]}\right)\int_0^{\infty}\exp\left(-\frac{y^2+2y\cdot\frac{1+2(\frac{\lambda}{\mu})^{\frac{1}{2}}}{(1+(\frac{\mu}{\lambda})^{\frac{1}{2}})^{\frac{1}{2}}}}{2}\right)dy\\
&\rightarrow&0\ \ {\rm as}\ \lambda\rightarrow\infty,
\end{eqnarray*}
and for $\lambda>\mu$,
\begin{eqnarray*}
&&\exp\left(\frac{2\lambda}{\mu}\right)\Phi\left(\frac{1-2(\frac{\lambda}{\mu})^{\frac{1}{2}}}{(1-(\frac{\mu}{\lambda})^{\frac{1}{2}})^{\frac{1}{2}}}\right)\\
&=&\frac{1}{\sqrt{2\pi}}\exp\left(\frac{2\lambda}{\mu}\right)\int_{\frac{-1+2(\frac{\lambda}{\mu})^{\frac{1}{2}}}{(1-(\frac{\mu}{\lambda})^{\frac{1}{2}})^{\frac{1}{2}}}}^{\infty}\exp\left(-\frac{y^2}{2}\right)dy\\
&=&\frac{1}{\sqrt{2\pi}}\exp\left(\frac{2\lambda}{\mu}\right)\int_{0}^{\infty}\exp\left(-\frac{\left(y+\frac{-1+2(\frac{\lambda}{\mu})^{\frac{1}{2}}}{(1-(\frac{\mu}{\lambda})^{\frac{1}{2}})^{\frac{1}{2}}}\right)^2}{2}\right)dy\\
&=&\frac{1}{\sqrt{2\pi}}\exp\left(-\frac{1}{2[1-(\frac{\mu}{\lambda})^{\frac{1}{2}}]}\right)\int_0^{\infty}\exp\left(-\frac{y^2+2y\cdot\frac{-1+2(\frac{\lambda}{\mu})^{\frac{1}{2}}}{(1-(\frac{\mu}{\lambda})^{\frac{1}{2}})^{\frac{1}{2}}}}{2}\right)dy\\
&\rightarrow&0\ \ {\rm as}\ \lambda\rightarrow\infty.
\end{eqnarray*}
Then,
\begin{eqnarray*}
&&\lim_{\lambda\rightarrow\infty}\Bigg[\Phi\left(\frac{1}{(1+(\frac{\mu}{\lambda})^{\frac{1}{2}})^{\frac{1}{2}}}\right)+\exp\left(\frac{2\lambda}{\mu}\right)\Phi\left(-\frac{1+2(\frac{\lambda}{\mu})^{\frac{1}{2}}}{(1+(\frac{\mu}{\lambda})^{\frac{1}{2}})^{\frac{1}{2}}}\right)\\
&&\ \ \ \ \ \ \ \  -\Phi\left(-\frac{1}{(1-(\frac{\mu}{\lambda})^{\frac{1}{2}})^{\frac{1}{2}}}\right)-\exp\left(\frac{2\lambda}{\mu}\right)\Phi\left(\frac{1-2(\frac{\lambda}{\mu})^{\frac{1}{2}}}{(1-(\frac{\mu}{\lambda})^{\frac{1}{2}})^{\frac{1}{2}}}\right)\Bigg]\\
&=&\Phi(1)-\Phi(-1)\\
&\approx&0.6827.
\end{eqnarray*}

Define
$$
y=\left(\frac{\mu}{\lambda}\right)^{\frac{1}{2}}.
$$
We have
\begin{eqnarray*}
&&\Phi\left(\frac{1}{(1+(\frac{\mu}{\lambda})^{\frac{1}{2}})^{\frac{1}{2}}}\right)+\exp\left(\frac{2\lambda}{\mu}\right)\Phi\left(-\frac{1+2(\frac{\lambda}{\mu})^{\frac{1}{2}}}{(1+(\frac{\mu}{\lambda})^{\frac{1}{2}})^{\frac{1}{2}}}\right)\\
&=&\Phi\left(\frac{1}{(1+y)^{\frac{1}{2}}}\right)+\exp\left(\frac{2}{y^2}\right)\Phi\left(-\frac{1+2y^{-1}}{(1+y)^{\frac{1}{2}}}\right)\\
&:=&J_1(y),
\end{eqnarray*}
and
\begin{eqnarray*}
&&\Phi\left(-\frac{1}{(1-(\frac{\mu}{\lambda})^{\frac{1}{2}})^{\frac{1}{2}}}\right)+\exp\left(\frac{2\lambda}{\mu}\right)\Phi\left(\frac{1-2(\frac{\lambda}{\mu})^{\frac{1}{2}}}{(1-(\frac{\mu}{\lambda})^{\frac{1}{2}})^{\frac{1}{2}}}\right)\\
&=&\Phi\left(-\frac{1}{(1-y)^{\frac{1}{2}}}\right)+\exp\left(\frac{2}{y^2}\right)\Phi\left(\frac{1-2y^{-1}}{(1-y)^{\frac{1}{2}}}\right)\\
&:=&J_2(y).
\end{eqnarray*}
To complete the proof, we need only show that $J_1(y)$ is an increasing function on $y\in(0,\infty)$ and $J_2(y)$ is a decreasing function on $y\in(0,1)$.

For $y>0$, we have
\begin{eqnarray*}
\frac{dJ_1}{dy}>0
&\Leftrightarrow&-\frac{1}{2\sqrt{2\pi}(1+y)^{\frac{3}{2}}}\exp\left(-\frac{1}{2(1+y)}\right)-\frac{4}{y^3}\exp\left(\frac{2}{y^2}\right)\Phi\left(-\frac{1+2y^{-1}}{(1+y)^{\frac{1}{2}}}\right)\\
&&+\frac{y^2+6y+4}{2\sqrt{2\pi}y^2(1+y)^{\frac{3}{2}}}\exp\left(-\frac{1}{2(1+y)}\right)>0\\
&\Leftrightarrow&\frac{3y+2}{\sqrt{2\pi}y^2(1+y)^{\frac{3}{2}}}\exp\left(-\frac{1}{2(1+y)}\right)-\frac{4}{y^3}\exp\left(\frac{2}{y^2}\right)\Phi\left(-\frac{1+2y^{-1}}{(1+y)^{\frac{1}{2}}}\right)>0\\
&\Leftrightarrow&\frac{(3y+2)y}{4\sqrt{2\pi}(1+y)^{\frac{3}{2}}}\exp\left(-\frac{(1+2y^{-1})^2}{2(1+y)}\right)>\Phi\left(-\frac{1+2y^{-1}}{(1+y)^{\frac{1}{2}}}\right),
\end{eqnarray*}
and for $0<y<1$, we have
\begin{eqnarray*}
\frac{dJ_2}{dy}<0
&\Leftrightarrow&\frac{1}{2\sqrt{2\pi}(1-y)^{\frac{3}{2}}}\exp\left(-\frac{1}{2(1-y)}\right)+\frac{4}{y^3}\exp\left(\frac{2}{y^2}\right)\Phi\left(\frac{1-2y^{-1}}{(1-y)^{\frac{1}{2}}}\right)\\
&&-\frac{y^2-6y+4}{2\sqrt{2\pi}y^2(1-y)^{\frac{3}{2}}}\exp\left(-\frac{1}{2(1-y)}\right)>0\\
&\Leftrightarrow&\frac{3y-2}{\sqrt{2\pi}y^2(1-y)^{\frac{3}{2}}}\exp\left(-\frac{1}{2(1-y)}\right)+\frac{4}{y^3}\exp\left(\frac{2}{y^2}\right)\Phi\left(\frac{1-2y^{-1}}{(1-y)^{\frac{1}{2}}}\right)>0\\
&\Leftrightarrow&\Phi\left(\frac{1-2y^{-1}}{(1-y)^{\frac{1}{2}}}\right)>\frac{(2-3y)y}{4\sqrt{2\pi}(1-y)^{\frac{3}{2}}}\exp\left(-\frac{(1-2y^{-1})^2}{2(1-y)}\right).
\end{eqnarray*}
For $y=1+r$ with $r\ge 0$, by H\"older's inequality, we get
\begin{eqnarray*}
&&\frac{(3y+2)y}{4\sqrt{2\pi}(1+y)^{\frac{3}{2}}}\exp\left(-\frac{(1+2y^{-1})^2}{2(1+y)}\right)>\Phi\left(-\frac{1+2y^{-1}}{(1+y)^{\frac{1}{2}}}\right)\nonumber\\
&\Leftrightarrow&\frac{(3y+2)y}{4(1+y)^{\frac{3}{2}}}>\int_0^{\infty}\exp\left(-\frac{\left(z+\frac{2(1+2y^{-1})}{(1+y)^{\frac{1}{2}}}\right)z}{2}\right)dz\nonumber\\
&\Leftarrow&\frac{(3y+2)y}{4(1+y)^{\frac{3}{2}}}>\left(\int_0^{\infty}\exp\left(-\frac{3z^2}{2}\right)dz\right)^{\frac{1}{3}}\left(\int_0^{\infty}\exp\left(-\frac{(1.5)(1+2y^{-1})z}{(1+y)^{\frac{1}{2}}}\right)dz\right)^{\frac{1}{1.5}}\nonumber\\
&\Leftrightarrow&\frac{(3y+2)y}{4(1+y)^{\frac{3}{2}}}>\frac{2^{\frac{1}{2}}\pi^{\frac{1}{6}}(1+y)^{\frac{1}{3}}}{3^{\frac{5}{6}}(1+2y^{-1})^{\frac{1}{1.5}}}\\
&\Leftrightarrow&3^5y^2(3y+2)^6(y+2)^4>2^{15}\pi(1+y)^{11}\\
&\Leftarrow&3^5y^2(3y+2)^6(y+2)^4>2^{17}(1+y)^{11}\\
&\Leftrightarrow&39111419 + 655929992 r + 2993031230 r^2 + 6991383720 r^3 +
 10103089845 r^4  \\
&&+ 9825833424 r^5+ 6700328484 r^6 + 3256962000 r^7 +
 1126305717 r^8\\
&& + 271238056 r^9 + 43292734 r^{10} + 4120456 r^{11} +
 177147 r^{12}>0.
\end{eqnarray*}
Hence, to complete the proof, we need only prove the following two inequalities:
\begin{eqnarray}\label{18aa}
\frac{(3y+2)y}{4\sqrt{2\pi}(1+y)^{\frac{3}{2}}}\exp\left(-\frac{(1+2y^{-1})^2}{2(1+y)}\right)>\Phi\left(-\frac{1+2y^{-1}}{(1+y)^{\frac{1}{2}}}\right),\ \ \ \ 0<y<1,
\end{eqnarray}
and
\begin{eqnarray}\label{18bb}
\Phi\left(\frac{1-2y^{-1}}{(1-y)^{\frac{1}{2}}}\right)>\frac{(2-3y)y}{4\sqrt{2\pi}(1-y)^{\frac{3}{2}}}\exp\left(-\frac{(1-2y^{-1})^2}{2(1-y)}\right),\ \ \ \ 0<y<\frac{2}{3}.
\end{eqnarray}

Denote the complementary error function by (cf. \cite{Wiki4})
$$
{\rm erfc}(x)=\frac{2}{\sqrt{\pi}}\int_x^{\infty}e^{-t^2}dt,\ \ \ \ x\in\mathbb{R}.
$$
We have
$$
\Phi(x)=\frac{1}{2}{\rm erfc}\left(-\frac{x}{\sqrt{2}}\right).
$$
By integration by parts, we get the following asymptotic expansion (cf. \cite{Water} and \cite{Wiki4}):
\begin{eqnarray*}
{\rm erfc}(x)=\frac{1}{\sqrt{\pi}}e^{-x^2}\left(\frac{1}{x}-\frac{1}{2x^3}+\frac{1\cdot 3}{2^2x^5}-\cdots+(-1)^{n-1}\frac{(2n-3)!!}{2^{n-1}x^{2n-1}}\right)+(-1)^n\frac{(2n-1)!!}{2^{n-1}\sqrt{\pi}}\int_x^{\infty}\frac{e^{-t^2}}{t^{2n}}dt.
\end{eqnarray*}
Then,
\begin{eqnarray}\label{18cc}
\Phi(-x)=\frac{1}{\sqrt{2\pi}}e^{-\frac{x^2}{2}}\left(\frac{1}{x}-\frac{1}{x^3}+\frac{1\cdot 3}{x^5}-\cdots+(-1)^{n-1}\frac{(2n-3)!!}{x^{2n-1}}\right)+(-1)^n\frac{(2n-1)!!}{2^n\sqrt{\pi}}\int_{\frac{x}{\sqrt{2}}}^{\infty}\frac{e^{-t^2}}{t^{2n}}dt.\nonumber\\
&&
\end{eqnarray}

For $0<y<1$, by (\ref{18cc}), we get
\begin{eqnarray*}
&&\frac{(3y+2)y}{4\sqrt{2\pi}(1+y)^{\frac{3}{2}}}\exp\left(-\frac{(1+2y^{-1})^2}{2(1+y)}\right)>\Phi\left(-\frac{1+2y^{-1}}{(1+y)^{\frac{1}{2}}}\right)\\
&\Leftarrow&\frac{(3y+2)y}{4(1+y)^{\frac{3}{2}}}>\frac{(1+y)^{\frac{1}{2}}}{1+2y^{-1}}-\left(\frac{(1+y)^{\frac{1}{2}}}{1+2y^{-1}}\right)^3+3\left(\frac{(1+y)^{\frac{1}{2}}}{1+2y^{-1}}\right)^5\\
&\Leftrightarrow&\frac{(3y+2)(y+2)}{4(1+y)^2}>1-\left(\frac{(1+y)^{\frac{1}{2}}}{1+2y^{-1}}\right)^2+3\left(\frac{(1+y)^{\frac{1}{2}}}{1+2y^{-1}}\right)^4\\
&\Leftrightarrow&\left(\frac{(1+y)^{\frac{1}{2}}}{1+2y^{-1}}\right)^2>\frac{y^2}{4(1+y)^2}+3\left(\frac{(1+y)^{\frac{1}{2}}}{1+2y^{-1}}\right)^4\\
&\Leftrightarrow&4(1+y)^3(y+2)^2>(y+2)^4+12y^2(1+y)^4\\
&\Leftrightarrow&y (32 + 64 y + 20 y^2 - 45 y^3 - 44 y^4 - 12 y^5)>0,
\end{eqnarray*}
which obviously holds for $0<y<1$. For $0<y<\frac{2}{3}$, by (\ref{18cc}), we get
\begin{eqnarray*}
&&\Phi\left(\frac{1-2y^{-1}}{(1-y)^{\frac{1}{2}}}\right)>\frac{(2-3y)y}{4\sqrt{2\pi}(1-y)^{\frac{3}{2}}}\exp\left(-\frac{(1-2y^{-1})^2}{2(1-y)}\right)\\
&\Leftarrow&\frac{(1-y)^{\frac{1}{2}}}{2y^{-1}-1}-\left(\frac{(1-y)^{\frac{1}{2}}}{2y^{-1}-1}\right)^3>\frac{(2-3y)y}{4(1-y)^{\frac{3}{2}}}\\
&\Leftrightarrow&4(1-y)^2(2-y)^2-4y^2(1-y)^3>(2-3y)(2-y)^3\\
&\Leftrightarrow&y^3 (8 - 11 y + 4 y^2)>0,
\end{eqnarray*}
which obviously holds for $0<y<\frac{2}{3}$. Therefore, inequalities (\ref{18aa}) and (\ref{18bb}) hold and the proof is complete.

\section{Numerical results for infinitely divisible discrete distributions and remarks}

\subsection{Geometric distribution}

Let $p>0$ and $X_{p}$ be a geometric random variable with parameter $p$. The probability mass function of $X_{p}$ is given by
$$
P\{X_{p}=k\}=p(1-p)^{k},\ \ \ \ k=0,1,2,\dots,
$$
and
$$
E[X_p]=\frac{1}{p}-1,\ \ \ \ {\rm Var}(X_p)=\frac{1-p}{p^2}.
$$
Define
$$
J_G(p):=P\left\{E[X_{p}]-\sqrt{{\rm Var}(X_{p})}<X_{p}\le E[X_{p}]+\sqrt{{\rm Var}(X_{p})}\right\}.
$$
Then, we have
$$
 P\left\{|X_{p}-E[X_{p}]|\le \sqrt{{\rm Var}(X_{p})}\right\}\ge J_G(p),\ \ \ \ p\in(0,1].
$$
By virtue of {\bf Mathematica}, we get
$$
 0.6827<\inf_p\{J_G(p)\}=\lim_{p\downarrow 0.75}J_G(p)=0.75<J_G(0.75)=0.9375.
$$
Below is the graph of the  function $J_G(p)$.
\begin{figure}[h]
\begin{center}
\scalebox{0.5}{\includegraphics{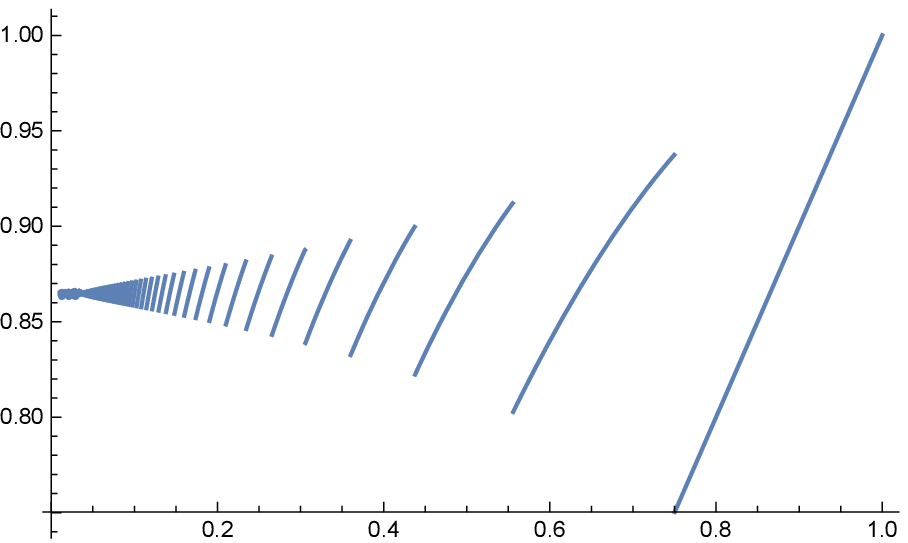}}
\end{center}
\end{figure}
\begin{center}
{\small Figure 3: Function $J_G(p)$.}
\end{center}

\subsection{Negative binomial distribution}

Let $n\ge2$, $p>0$ and $X_{n,p}$ be a negative binomial random variable with parameters $n$ and $p$. In general, we do not have
$$
I_{NB}(n,p):= P\left\{|X_{n,p}-E[X_{n,p}]|\le \sqrt{{\rm Var}(X_{n,p})}\right\}\ge 0.6827.
$$
For example,
$$
I_{NB}(2,0.45)=0.6339326< 0.6827.
$$
Below is the  graph of the function $I_{NB}(2,p)$.
\newpage\begin{figure}[h]
\begin{center}
\scalebox{0.5}{\includegraphics{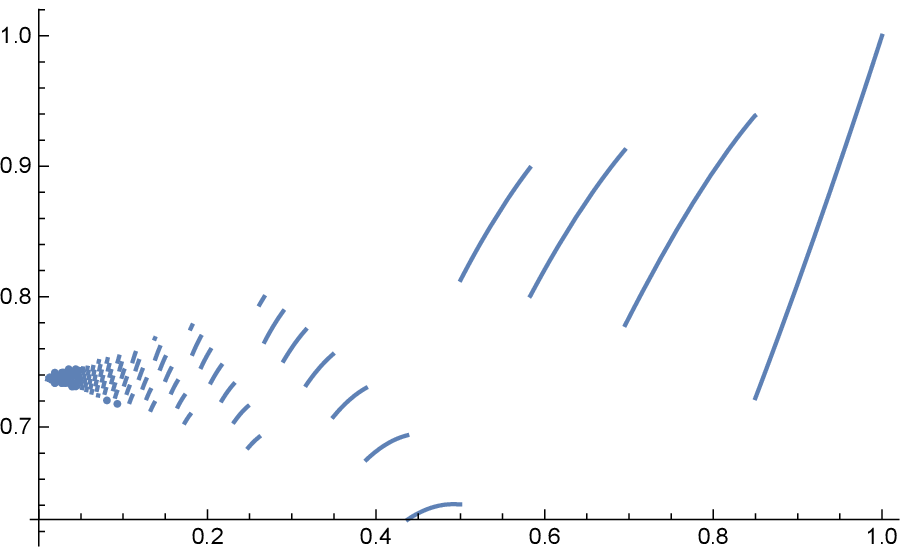}}
\end{center}
\end{figure}
\begin{center}
{\small Figure 4: Function $I_{NB}(2,p)$.}
\end{center}

Motivated by the above example, for the negative binomial distribution, we suggest consider the following slightly modified inequality:
$$
J_{NB}(n,p):=P\left\{\lfloor E[X_{n,p}]-\sqrt{{\rm Var}(X_{n,p})}\rfloor\le X_{n,p}\le E[X_{n,p}]+\sqrt{{\rm Var}(X_{n,p})}\right\}> 0.6827.
$$
Hereafter $\lfloor x\rfloor$ denotes the greatest integer less than or equal to $x$.
By virtue of {\bf Mathematica}, we get
$$
J_{NB}(n,p)>0.6827,\ \ \ \ p\in(0,1],\ n\ge2.
$$
Below are graphs of the  function $J_{NB}(n,p)$ for $n=2,3,10,1000$.
\begin{figure}[h]
\begin{center}
\scalebox{0.5}{\includegraphics{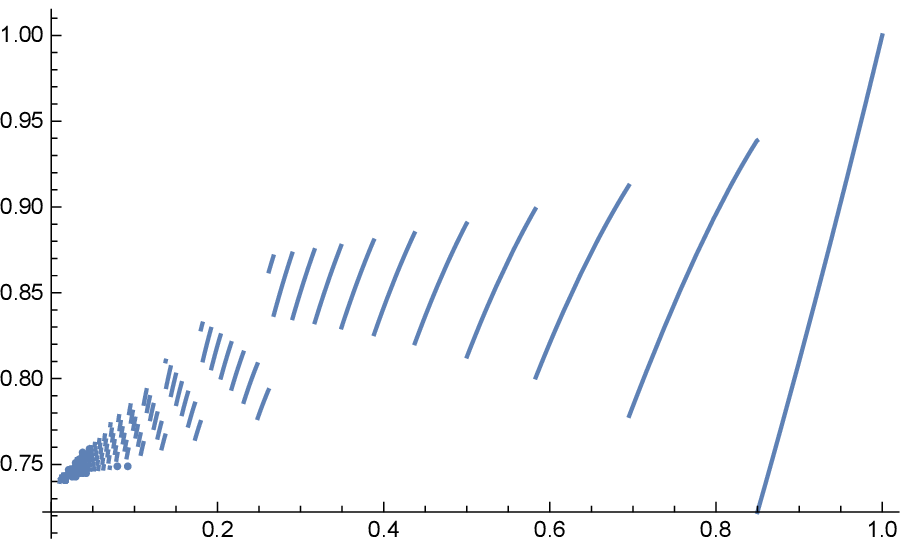}}
\end{center}
\end{figure}
\begin{center}
{\small Figure 5: Function $J_{NB}(2,p)$.}
\end{center}

\begin{figure}[h]
\begin{center}
\scalebox{0.5}{\includegraphics{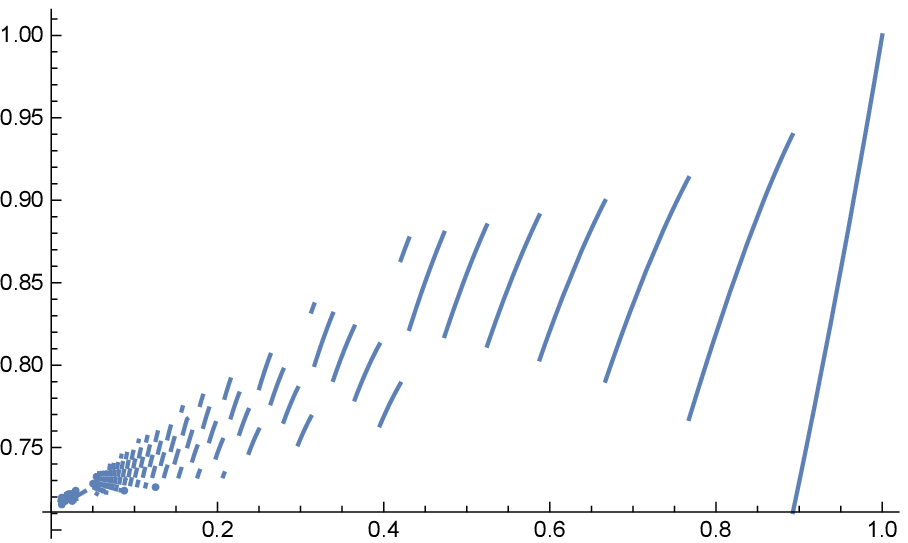}}
\end{center}
\end{figure}
\begin{center}
{\small Figure 6: Function $J_{NB}(3,p)$.}
\end{center}

\newpage\begin{figure}[h]
\begin{center}
\scalebox{0.5}{\includegraphics{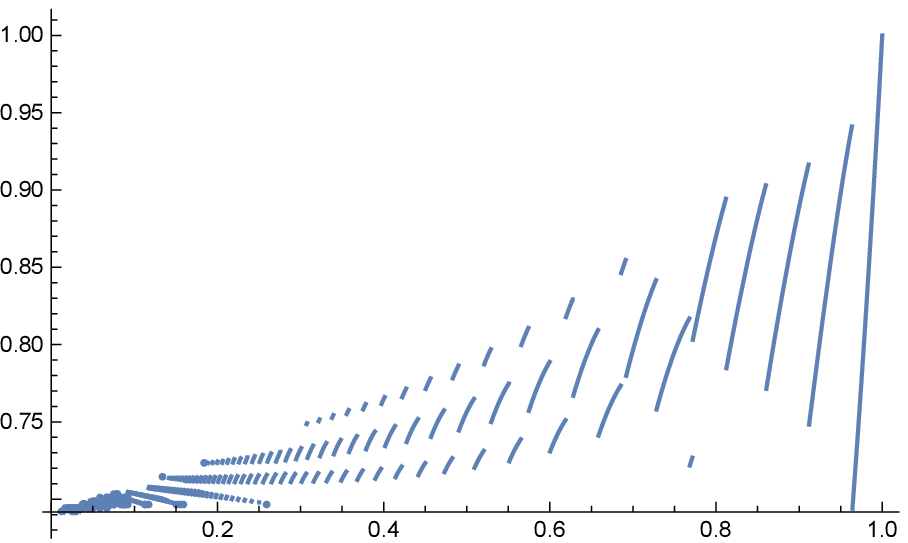}}
\end{center}
\end{figure}
\begin{center}
{\small Figure 7: Function $J_{NB}(10,p)$.}
\end{center}

\begin{figure}[h]
\begin{center}
\scalebox{0.5}{\includegraphics{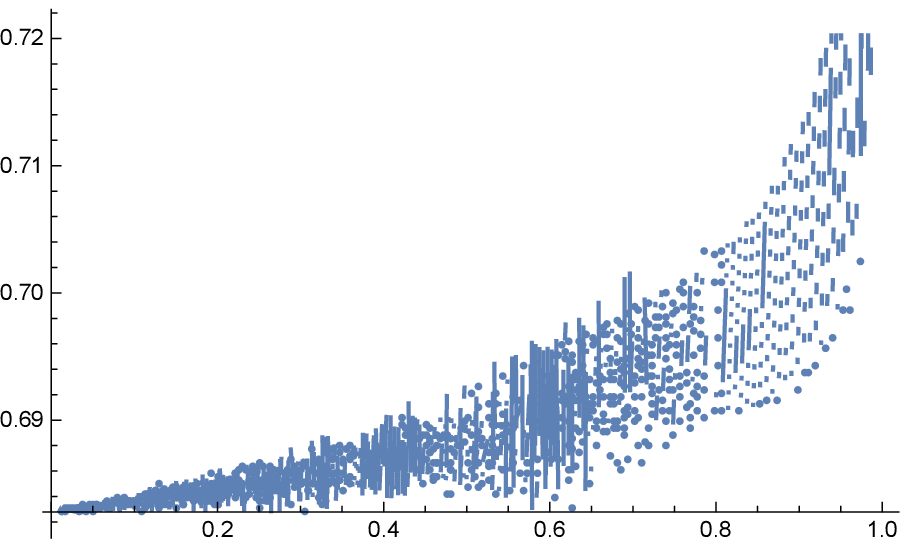}}
\end{center}
\end{figure}
\begin{center}
{\small Figure 8: Function $J_{NB}(1000,p)$.}
\end{center}

\subsection{Poisson distribution}\label{Poisson1}

Let $\lambda>0$ and $X_{\lambda}$ be a Poisson random variable with parameter $\lambda$.  In general, we do not have
$$
I_{P}(\lambda):= P\left\{|X_{\lambda}-E[X_{\lambda}]|\le \sqrt{{\rm Var}(X_{\lambda})}\right\}\ge 0.6827.
$$
For example,
$$
I_{P}(3)= 0.616115< 0.6827.
$$

Motivated by the above example, for the Poisson distribution, we suggest consider the following slightly modified inequality:
$$
J_{P}(\lambda):=P\left\{\lfloor E[X_{\lambda}]-\sqrt{{\rm Var}(X_{\lambda})}\rfloor\le X_{\lambda}\le \lceil E[X_{\lambda}]+\sqrt{{\rm Var}(X_{\lambda})}\rceil \right\}> 0.6827.
$$
Hereafter $\lceil x\rceil$ denotes the smallest integer greater than or equal to $x$.
By virtue of {\bf Mathematica}, we get
$$
J_{P}(\lambda)>0.6827,\ \ \ \ \lambda>0.
$$
Below is the graph of the  function $J_P(\lambda)$ for $\lambda\in(0,100]$.
\newpage\begin{figure}[h]
\begin{center}
\scalebox{0.5}{\includegraphics{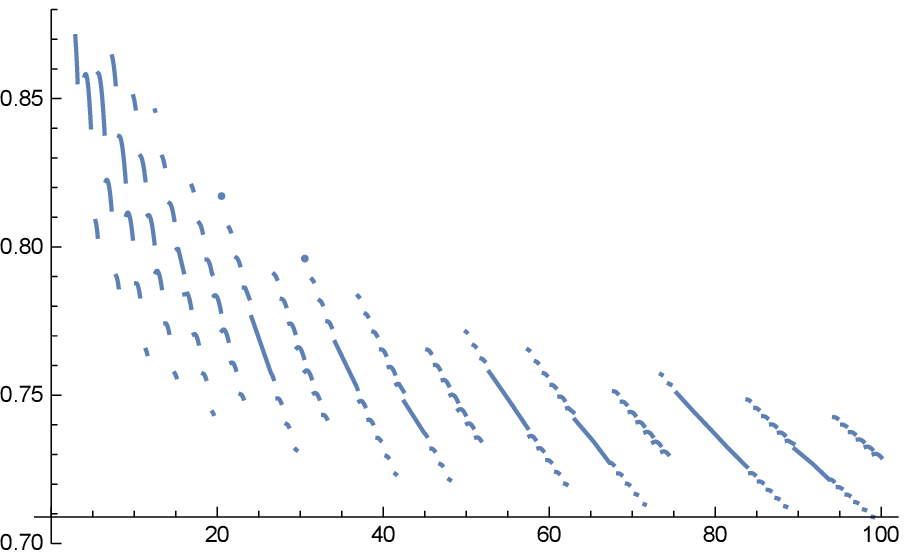}}
\end{center}
\end{figure}
\begin{center}
{\small Figure 9: Function $J_P(\lambda)$.}
\end{center}

\subsection{Remarks}

Let $X_{3}$ be a Poisson random variable with parameter $3$. By \S \ref{Poisson1}, we know that
$$
P\left\{|X_3-E[X_{3}]|\le \sqrt{{\rm Var}(X_{3})}\right\}=0.616115< 0.6827.
$$
Let $B$ be a standard normal random variable which is independent of $X_{3}$. Define
$$
X_{B,\varepsilon}:=\varepsilon B+X_{3}.
$$
Then, we  have
$$
\lim_{\varepsilon\rightarrow0}P\left\{|X_{B,\varepsilon}-E[X_{B,\varepsilon}]|\le \sqrt{{\rm Var}(X_{B,\varepsilon})}\right\}=0.616115< 0.6827.
$$
Hence, inequality (\ref{inf300}) does not hold for all infinitely divisible continuous distributions.

For $n\in\mathbb{N}$, define
$$
\nu_{n}(dx)=\frac{3n}{2}\cdot1_{\left[1-\frac{1}{n},1+\frac{1}{n}\right]}(x)dx.
$$
Let $Y_n$ be a compound Poisson random variable with L\'evy measure $\nu_{n}$. Then, $Y_n$ converges to $X_{3}$ in distribution as $n\rightarrow\infty$. Hence, $Y_n$ does not satisfy (\ref{inf300}) at least if $n$ is large enough. This simple example shows that inequality  (\ref{inf300}) might not hold if the L\'evy measure of the infinitely divisible random variable is finite, even if this measure is absolutely continuous with respect to the Lebesgue measure.

It deserves considering under what conditions an  infinitely divisible random variable with infinite L\'evy measure satisfies (\ref{inf300}). Also, it is interesting to investigate this inequality  for general (not necessarily parametric) infinitely divisible continuous distributions and consider, to what extent, it can be established for non-infinitely divisible continuous distributions.
\vskip 0.5cm
{ \noindent {\bf\large Acknowledgements}\quad  This work was supported by the National Natural Science Foundation of China (No. 12171335), the Science Development Project of Sichuan University (No. 2020SCUNL201) and the Natural Sciences and Engineering Research Council of Canada (No. 4394-2018).

\end{document}